\documentclass[12pt]{article}
\pdfoutput=1
\usepackage{jhep-mod}
\usepackage{bm}
\usepackage{amssymb}
\usepackage{pifont}
\usepackage{epstopdf} 
\usepackage[applemac]{inputenc}
\title{Primes and the Lambert $W$ function}
\author{Matt Visser}
\affiliation{School of Mathematics and Statistics, \\
Victoria University of Wellington, PO Box 600, Wellington 6140, New Zealand}
\emailAdd{matt.visser@sms.vuw.ac.nz}
\abstract{

\noindent
The Lambert $W$ function, implicitly defined by $W(x)\, e^{W(x)}=x$, is a ``new'' special function that has recently been the subject of an extended upsurge in interest and applications. In this note, I point out that the Lambert $W$ function can also be used to gain a new perspective on the distribution of primes. 

\bigskip
\noindent
1 February 2017; \LaTeX-ed \today
}
\keywords{Lambert $W$ function; prime counting function; the $n$'th prime.\\[2pt]
arXiv: \qquad\;\;\,1311.2324 [math.NT].
} 

\begin{document}
\maketitle
\def\R{{\mathbb{R}}}
\newtheorem{theorem}{Theorem}
\newtheorem{corollary}{Corollary}
\newtheorem{lemma}{Lemma}
\section{Introduction}
\label{S:intro}

The Lambert $W$ function, implicitly defined by the relation $W(x)\, e^{W(x)}=x$, has a long and quite convoluted 250-year history, but only recently has it become common to view this particular function as one of the standard special functions of mathematics~\cite{Corless}. Applications range widely~\cite{Corless, Valluri:00}, from combinatorics (for instance in the enumeration of rooted trees)~\cite{Corless},  to delay differential equations~\cite{Corless}, to falling objects subject to linear drag~\cite{Vial:12}, to the evaluation of the numerical constant in Wien's displacement law~\cite{Stewart:11, Stewart:12}, to quantum statistics~\cite{Valluri:09}, to constructing the ``tortoise'' coordinate for Schwarzschild black holes~\cite{tortoise}, \emph{etcetera}. In this brief note I will indicate some apparently new applications of the Lambert $W$ function to the distribution of primes, specifically to the prime counting function $\pi(x)$ and estimating the $n$'th prime $p_n$.

\section{The prime counting function $\pi(x)$}\label{S:Pi}
\subsection{Upper bound\label{SS:Pi-U}}
\begin{theorem}
The prime counting function $\pi(x)$ satisfies
\begin{equation}
\pi(x) < {x\over W(x)} =  e^{W(x)};   \qquad (\forall x \geq 0).
\end{equation}
\end{theorem}
\paragraph{Proof:} First note that $ x \geq p_{\pi(x)}$. Second recall the standard result that $p_n > n\ln n$ for $n \geq 1$. (See Rosser~\cite{Rosser:38}, or any standard reference book on  prime numbers, for example~\cite{Ribenboim:91, Ribenboim:96}.) Then
$x \geq p_{\pi(x)} > \pi(x) \ln\pi(x)$, so we have $x > \pi(x) \ln\pi(x)$. 
Invert, noting that the RHS is monotone increasing, to see that $\pi(x) < x/ W(x)$, certainly for $\pi(x)\geq 1$ (corresponding to $x\geq2$). Then explicitly check validity of the inequality on the domain $x\in[0,2)$. Finally, use the definition of the Lambert $W$ function to note $ x/ W(x) = e^{W(x)}$. \hfill{$\Box$}

\subsection{Asymptotics\label{SS:Pi-A}}
\begin{theorem}
The prime number theorem, $\pi(x) \sim x/\ln x$, is equivalent to the statement
\begin{equation}
\pi(x) \sim {x\over W(x)} =  e^{W(x)};   \qquad (x\to\infty).
\end{equation}
\end{theorem}
\paragraph{Proof:} Trivial. Note that asymptotically $W(x) \sim \ln x$. (The only potential subtlety is that we are using the principal branch of the $W$ function,  denoted $W_0(x)$ whenever there is any risk of confusion~\cite{Corless}.)
\hfill{$\Box$}

\noindent
So we have derived \emph{both} a strict upper bound on $\pi(x)$ \emph{and} an asymptotic equality. 

\subsection{Lower bound\label{SS:Pi-L}}

When it comes to developing an analogous lower bound on $\pi(x)$,  the situation is considerably more subtle. First consider the well-known upper bound on $p_n$~\cite{Rosser:41, Ribenboim:91, Ribenboim:96}: 
\begin{equation}
p_n < n(\ln n+ \ln\ln n) = n \ln( n \ln n );  \qquad (n\geq 6).
\label{E:lnln}
\end{equation}
Note that $x < p_{\pi(x)+1}$. 
Now consider the elementary inequality
\begin{equation}
\ln x \leq {x\over e},
\end{equation}
with equality only at $x=e$,
and observe that this implies
\begin{equation}
\ln x = {\ln(x^\epsilon)\over\epsilon} \leq  {x^\epsilon\over\epsilon\; e};  \qquad (\forall \epsilon > 0),
\end{equation}
now with equality only at $x= e^{1/\epsilon}$.
This inequality explicitly captures the well-known fact that the logarithm grows less rapidly than any positive power. Applied to the $n$'th prime this now implies
\begin{equation}
p_n < n \ln( n \ln n ) \leq  n \ln\left(n^{1+\epsilon}\over\epsilon\; e\right) = n \{(1+\epsilon)\ln n -1-\ln\epsilon \};  \qquad (n\geq 6; \;\;\forall \epsilon>0).
\end{equation}
This inequality is much weaker than equation (\ref{E:lnln}), but much more tractable. 
Using logic identical to that of Theorem 1, it is easy to convert this into the inequality
\begin{equation}
x\leq p_{\pi(x)+1} <  [\pi(x)+1] \ln\left\{ [\pi(x)+1]^{1+\epsilon}\over \epsilon\; e \right\};  \qquad (\pi(x) \geq 5; \;\;\forall \epsilon>0).
\end{equation}
This is now easily inverted to obtain:
\begin{theorem}
The prime counting function $\pi(x)$ satisfies \emph{($x \geq 11; \;\;\forall \epsilon>0$):}
\begin{equation}
\pi(x) >  { \displaystyle {x\over 1+\epsilon} 
\over 
W\left( \displaystyle{ x\over 1+\epsilon} \; (\epsilon\; e)^{-1/(1+\epsilon)} \right)} - 1
\,=\, (\epsilon\; e)^{1/(1+\epsilon)}  \exp  W\left( \displaystyle{ x\over 1+\epsilon} \; (\epsilon\; e)^{-1/(1+\epsilon)} \right).
\end{equation}
\end{theorem}
For any fixed $\epsilon>0$,  this inequality holds at least for $\pi(x)\geq 5$, corresponding to $x\geq 11$. But, depending on the specific value $\epsilon$, the domain of validity may actually be larger.  That is,  $\forall\epsilon>0$ the inequality holds for $\pi(x)\geq n_0(\epsilon)$ with $n_0(\epsilon)\leq 5$, corresponding to $x\geq x_0(\epsilon)$ with $x_0(\epsilon)<11$.
For instance, numerically solving ${1\over 1+\epsilon} \; (\epsilon\; e)^{-1/(1+\epsilon)} =1$ gives $\epsilon_* = 0.2711715619...$, whence ${1\over 1+\epsilon_*} = (\epsilon\; e)^{1/(1+\epsilon_*)} = 0.7866758744...$. 
\begin{corollary}
After explicitly checking the domain of validity we have:
\begin{equation}
  \pi(x) >  { x \over (1+\epsilon_*) \; W\left( x \right)} - 1 \,=\,
  {\exp W\left(x \right)\over 1 +\epsilon_*}-1; \qquad (x \geq 3).
\end{equation}
Where numerically ${1\over 1+\epsilon_*} = 0.7866758744...$.
\end{corollary} 
Other corollaries may be determined analytically:
\begin{corollary}
Setting $\epsilon=1$, and explicitly checking the domain of validity, we have:
\begin{equation}
  \pi(x) >  { \displaystyle {x/2} \over W\left( \displaystyle{ x \over2 e^{1/2}} \;  \right)} - 1 \,=\,
  e^{1/2} \, \exp W\left( \displaystyle{ x \over 2e^{1/2}} \;  \right)-1; \qquad (x \geq 6).
\end{equation}
\end{corollary} 
\begin{corollary}
Choosing the specific case $\epsilon = e^{-1}$, and explicitly checking the domain of validity,  we have:
\begin{equation}
\pi(x) >  { \displaystyle {x\over 1+e^{-1}} \over W\left( \displaystyle{ x\over 1+e^{-1}} \; \right)} - 1
\,=\,
\exp W\left( \displaystyle{ x\over 1+e^{-1}} \; \right) - 1; \qquad (x\geq 5).
\end{equation}
\end{corollary}
\begin{corollary}
Choosing the specific case  $\epsilon = e^{-3}$,  the domain of validity is the entire positive half line  \emph{($x\geq 0$):}\begin{equation}
\pi(x) >  { \displaystyle {x\over 1+e^{-3}} \over W\left( \displaystyle{ x\over 1+e^{-3}} \; e^{2/(1+e^{-3})}\right)} - 1
\,=\,
e^{-2/(1+e^{-3})} \,\exp W\left( \displaystyle{ x\over 1+e^{-3}} \; e^{2/(1+e^{-3})}\right)- 1. 
\end{equation}
\end{corollary}
This fourth corollary, ($\epsilon=e^{-3}$), exhibits somewhat poorer bounding performance at intermediate values of $x$, but eventually  overtakes the third corollary, ($\epsilon=e^{-1}$), once $x\approx e^{2e+3}\approx 4600$, and then asymptotically provides a better bound. 
Numerous variations on this theme can also be constructed, amounting to different ways of approximating the logarithms appearing in equation (\ref{E:lnln}).

In summary, we have used the Lambert $W$ function to obtain a number of bounds, and some general classes of bounds, on the prime counting function $\pi(x)$ in terms of the Lambert $W$ function $W(x)$. We shall now turn attention to the $n$'th prime $p_n$. 

\section{The $n$'th prime}\label{S:p_n}
\subsection{Upper bound\label{SS:p_n-L}}

\begin{theorem}
The  $n$'th prime $p_n$ satisfies
\begin{equation}
p_n < - n \;W_{-1}\left(-{1\over n}\right) = n \left| W_{-1}\left(-{1\over n}\right)\right| \qquad (n \geq 4).
\end{equation}
Here $ W_{-1}(x)$ is the second real branch of the Lambert $W$ function, defined on the domain $x\in [-1/e, 0)$. 
\end{theorem}
\paragraph{Proof:}
We start from the fact that
$n \geq p_n /\ln p_n$,  this inequality certainly being valid for $p_n\geq 7$, corresponding to $n\geq 4$~\cite{Dusart-1999}. Inverting, (and appealing to the monotonicity of $x/\ln x$), we have $p_n < - n \;W_{-1}(-1/n)$, certainly for $n \geq 4$. Explicitly inspecting $n\in\{1,2,3\}$ shows that the actual range of validity is indeed $n\geq4$.   \hfill{$\Box$}

\begin{corollary}
\begin{equation}
p_n < - n \;W_{-1}\left(-{1\over n+e}\right) =  n \left|W_{-1}\left(-{1\over n+e}\right)\right| \qquad (n \geq 1).
\end{equation}
\end{corollary}
\paragraph{Proof:} Note that $-W_{-1}(x)$ is monotone increasing on $[-e^{-1},0)$. 
So we see that $-W_{-1}(-1/[n+e]) > -W_{-1}(-1/n)$, and the claimed inequality certainly holds for $n\geq 4$. For $n\in\{1,2,3\}$ verify the claimed inequality by explicit computation. \hfill{$\Box$}

\bigskip
\noindent
The virtue of this specific corollary is that it now holds for all positive integers. \\
There are many other variations on this theme that one could construct.


\subsection{Asymptotics\label{SS:p_n-A}}

\begin{theorem}
The prime number theorem, which can be written in the form $p_n \sim n\ln n$, is equivalent to the statement
\begin{equation}
p_n \sim - n \;W_{-1}\left(-{1\over n}\right) = n \left| W_{-1}\left(-{1\over n}\right)\right| ;   \qquad (n\to\infty).
\end{equation}
\end{theorem}
\paragraph{Proof:} Trivial. Consider the asymptotic result 
\begin{equation}
W_{-1}(x) =  \ln(-x) - \ln(-\ln(-x)) + o(1)  \qquad (x\to 0^-).
\end{equation}
Then
\begin{equation}
- n \;W_{-1}\left(-{1\over n}\right) = n \{ \ln n + \ln\ln n + o(1)\}.
\end{equation}
\paragraph{Comment:}
Note that use of the Lambert $W$ function automatically yields the first two terms of the Cesaro--Cippola asymptotic expansion~\cite{Cesaro, Cippola}:
\begin{equation}
p_n = n \{ \ln n + \ln\ln n -1 + o(1)\}.
\end{equation}
We can even obtain the first \emph{three} terms of the Cesaro--Cippola asymptotic expansion by refining the prime number theorem slightly as follows:
\begin{theorem}
\begin{equation}
p_n \sim - n \;W_{-1}\left(-{e\over n}\right)  = n \left| W_{-1}\left(-{e\over n}\right)\right| ;   \qquad (n\to\infty).
\end{equation}
\end{theorem}
Finally, we note that use of the Lambert $W$ function yields \emph{both} a strict upper bound \emph{and} an asymptotic result.

\subsection{Lower bound\label{SS:p_n-U}}

To obtain an lower bound on $p_n$ we start with an upper bound on $\pi(x)$. Consider for instance the standard result~\cite{Rosser:62}:
\begin{equation}
\pi(x) < {x\over\ln x - {3\over2}}; \qquad (x>e^{3/2})
\end{equation}
Note that the RHS of this inequality is monotone increasing for $x>e^{5/2}$. 
Now we always have $p_{\pi(x)} \leq x<p_{\pi(x)+1}$, so 
\begin{equation}
n < {p_{n+1} \over \ln p_{n+1} - {3\over2}}.
\end{equation}
This holds at the very least for $p_n > e^{5/2}$, corresponding to $n \geq 6$, but an explicit check shows that it in fact holds for $n\geq 2$.
This is perhaps more clearly expressed as
\begin{equation}
n -1 < {p_{n} \over \ln p_{n} - {3\over2}};  \qquad (n \geq 3).
\end{equation}
Inverting, and noting the constraint arising from the domain of definition of $W_{-1}$, we now obtain:
\begin{theorem}
\begin{equation}
p_n > - (n-1) \; W_{-1}\left(-{e^{3/2}\over n-1}\right) =  (n-1) \left| W_{-1}\left(-{e^{3/2}\over n-1}\right)\right| ; \qquad ( n \geq 14).
\end{equation}
\end{theorem}
There are many other variations on this theme that one could construct.

\subsection{Implicit bounds\label{SS:p_n-implicit}}

It is an old result (see for example Rosser~\cite{Rosser:41}) that $\forall \epsilon >0, \; \exists N(\epsilon): \forall n \geq N(\epsilon)$
\begin{equation}
{x\over\ln x - 1 + \epsilon } <\pi(x) < {x\over\ln x - 1 - \epsilon }.
\end{equation}
Without an explicit calculation of $N(\epsilon)$ these bounds are qualitative, rather than quantitative. Nevertheless it may be of interest to point out that a minor variant of the arguments above immediately yields:
\begin{theorem}
$\forall \epsilon >0, \; \exists M(\epsilon): \forall n \geq M(\epsilon)$
\begin{equation}
- n \; W_{-1}\left(-{e^{1-\epsilon}\over n}\right) > p_n > - (n-1) \; W_{-1}\left(-{e^{1+\epsilon}\over n-1}\right).
\end{equation}
\end{theorem}
It is now ``merely'' a case of estimating $M(\epsilon)$ to turn these into explicit bounds. 
We have already seen that $\epsilon =1$ provides a widely applicable upper bound, and $\epsilon=1/2$ a widely applicable lower bound.
Taking $\epsilon\to 0$ now makes it clear why
\begin{equation}
p_n \sim - n \;W_{-1}\left(-{e\over n}\right) ;   \qquad (n\to\infty),
\end{equation}
is such a good asymptotic estimate for $p_n$. 

\section{Discussion}\label{S:Discussion}

While the calculations carried out above are very straightforward, almost trivial, it is perhaps the shift of viewpoint that is more interesting. The Lambert $W$ function provides (in this context) a ``new'' special function to work with, one which may serve to perhaps simplify and unify many otherwise disparate results. It is perhaps worth noting that the infamous ``$\ln\ln x$'' terms that infest the analytic theory of prime numbers will automatically appear as the sub-leading terms in asymptotic expansions of the Lambert $W$ function. 

\section*{Acknowledgments}

This research was supported by the Marsden Fund, and by a James Cook fellowship, both administered by the Royal Society of New Zealand.  

\appendix
\renewcommand{\theequation}{A.\arabic{equation}}  
\setcounter{equation}{0}  
\section{Appendix}\label{S:Appendix}

The Lambert $W$ function is a multi-valued complex function defined implicitly by~\cite{Corless} 
\begin{equation}
W(x)\; e^{W(x)}=x,
\end{equation}
There are two real branches: $W_0(x)$ defined for $x\in[-e^{-1},\infty)$, and $W_{-1}(x)$ defined for $x\in[e^{-1},0)$. These two branches meet at the common point $W_0(-e^{-1})=W_{-1}(-e^{-1})=-1$.  It is common to use $W(x)$ in place of $W_0(x)$ when there is no risk of confusion.

\noindent
Asymptotic expansions are:
\begin{equation}
W_0(x) =  \ln x - \ln\ln x + o(1); \qquad (x\to\infty);
\end{equation}
\begin{equation}
W_{-1}(x) =  \ln(-x) - \ln(-\ln(-x)) + o(1);  \qquad (x\to 0^-).
\end{equation}
More details, and a Taylor expansion for $|x|<e^{-1}$, can be found in Corless~\emph{et al,}~\cite{Corless}.
\enlargethispage{20pt}

\noindent
A key identity is:
\begin{equation}
\ln(a+bx)+cx=\ln d \qquad\implies\qquad x = {1\over c} \; W\left( {cd\over b} \exp\left[{ac\over b}\right]\right) - {a\over b}.
\end{equation}

\clearpage



\begin{thebibliography}{99}

\bibitem{Corless}
Corless, R. M.; Gonnet, G. H.; Hare, D. E. G.; Jeffrey, D. J.; and Knuth, D. E.,\\
``On the Lambert $W$ function", \\
Advances in Computational Mathematics {\bf5} (1996) 329--359. 
doi:10.1007/BF02124750.

\bibitem{Valluri:00}
S. R. Valluri, D. J. Jeffrey, R. M. Corless, \\
``Some Applications of the Lambert $W$ Function to Physics'',\\
Canadian Journal of Physics {\bf78} (2000) 823--831, 
doi:10.1139/p00-065

\bibitem{Vial:12}
Alexandre Vial, ``Fall with linear drag and Wien's displacement law: approximate solution and Lambert function'',
European Journal of Physics, {\bf33} (2012) 751,  doi:10.1088/0143-0807/33/4/751

\bibitem{Stewart:11}
Se\'an M Stewart, ``Wien peaks and the Lambert $W$ function'', \\
Revista Brasileira de Ensino de F\'isica {\bf33} (2011) 3308, {\sf www.sbfisica.org.br}

\bibitem{Stewart:12}
Se\'an M. Stewart,  ``Spectral Peaks and Wien's Displacement Law'', \\
Journal of Thermophysics and Heat Transfer, {\bf26} (2012) 689--692.
doi: 10.2514/1.T3789


\bibitem{Valluri:09}
S. R. Valluri, M. Gil, D. J. Jeffrey, and Shantanu Basu,\\
``The Lambert $W$ function and quantum statistics'',\\
Journal of Mathematical Physics {\bf50} (2009) 102103, doi:10.1063/1.3230482

\bibitem{tortoise}
Petarpa~Boonserm and Matt~Visser,\\
  ``Bounding the greybody factors for Schwarzschild black holes'',\\
  Phys.\ Rev.\ D {\bf 78} (2008) 101502
  [arXiv:0806.2209 [gr-qc]].
  
\bibitem{Rosser:38}
J. Barkley Rosser, ``The $n$'th Prime is Greater than $n\ln n$", \\
Proc. London Math. Soc. {\bf45} (1938) 21--44.

\bibitem{Ribenboim:91}
Ribenboim, P. \emph{The Little Book of Big Primes}, (Springer-Verlag, New York, 1991).

\bibitem{Ribenboim:96}
Ribenboim, P. \emph{The New Book of Prime Number Records},\\
(Springer-Verlag, New York, 1996).

\bibitem{Rosser:41}
J. Barkley Rosser, ``Explicit Bounds for some functions of prime numbers'', \\
American  Journal of  Mathematics, {\bf 63}  (1941) 211--232.


\bibitem{Dusart-1999}
Pierre Dusart,
``The $k^{th}$ prime is greater than $k(\ln k +\ln\ln k -1)$ for $k\geq 2$'',\\
Mathematics of Computation {\bf 68} (1999) 411--415.

\bibitem{Cesaro}
Ernest Ces\`aro. ``Sur une formule empirique de M. Pervouchine", \\
Comptes rendus  {\bf119} (1894) 848--849. (French)
\enlargethispage{35pt}

\bibitem{Cippola}
M. Cipolla, ``La determinazione assintotica dell'$n^{imo}$ numero primo'', \\
Matematiche Napoli, {\bf3} (1902) 132--166. (Italian)

\bibitem{Rosser:62}
J. Barkley Rosser and Lowell Schoenfeld, \\
``Approximate Formulas for Some Functions of Prime Numbers", \\
Illinois J. Math. {\bf6} (1962) 64--97.

\end{thebibliography}
\end{document}